\documentclass[12pt]{article}
\usepackage{geometry}                
\usepackage{graphicx,amssymb}
\usepackage{graphicx}
\usepackage{rotating}
\usepackage{dcolumn}
\usepackage{amsfonts}
\usepackage{amssymb}
\usepackage{amsmath}
\usepackage{latexsym}
\usepackage{epsfig}
\usepackage{dsfont,eucal,mathrsfs}
\usepackage{hyperref}
\usepackage{verbatim}
\usepackage[textsize=footnotesize]{todonotes}
\usepackage{color}
\usepackage{amsfonts}
\usepackage{amsmath}
\usepackage{amssymb}
\usepackage{amsthm}
\usepackage{bm}

\newcommand{\keywords}[1]{\par\addvspace\baselineskip\noindent\textbf{Keywords:}\enspace\ignorespaces#1}
\newcommand{\AMSclassification}[1]{\par\addvspace\baselineskip\noindent\textbf{2010 Mathematics Subject Classification:}\enspace\ignorespaces#1}

\newtheorem{theorem}{Theorem}
\newtheorem{lemma}[theorem]{Lemma}
\newtheorem{propo}[theorem]{Proposition}
\newtheorem{coro}[theorem]{Corollary}
\newtheorem{definition}{Definition}

\begin{document}

\title{Bifurcations of mutually coupled equations in random graphs}
\author{
Eduardo Garibaldi\thanks{Supported by BREUDS and FAPESP 2012/18780-0.} \\UNICAMP -- Department of Mathematics \\ 13083-859 Campinas -- SP, Brazil \\ garibaldi@ime.unicamp.br \\~\\
Tiago Pereira\thanks{Supported by FP7 Marie Currie Actions (Project 303180) and IRSES--DynEurBraz.} \\Department of Mathematics \\ Imperial College London, South Kensington, London SW7 2BA  \\
tiago.pereira@imperial.ac.uk
}

\date{}

\maketitle

\begin{abstract}
We study the behavior of solutions of mutually coupled equations in heterogeneous random graphs.  Heterogeneity means that some equations receive many inputs whereas most of the equations are given only with a few connections.  Starting from a situation where the isolated equations are unstable, we prove that a heterogeneous interaction structure leads to the appearance of stable subspaces of solutions.   Moreover, we show that, for certain classes of heterogeneous networks, increasing the strength of interaction leads to a cascade of bifurcations in which the dimension of the stable subspace of solutions increases. We explicitly determine the bifurcation scenario in terms of the graph structure. 

\keywords{bifurcation, coupled equations, dichotomies, random graphs}

\AMSclassification{05C80, 34C15, 34F05, 34F10, 37C10, 60B20}

\end{abstract}

\section{Introduction}

The last decade has witnessed rapidly growing interest in dynamics of coupled dynamical systems  \cite{Barrat,Book,Young}.  In most applications, the interaction structure among elements is intricate \cite{Newman} and modeled by random graphs \cite{Bollobas,Chung}.  Empirical studies indicate that this interaction structure can have dramatic influences on the dynamical properties and the functioning of such systems \cite{Bullmore,PRL}. 

Recent studies show that disparate real-world networks display a heterogeneous connectivity -- while some nodes, called hubs, receive many connections, most of the nodes are poorly connected  \cite{Barrat,Bullmore}.  Such a connectivity structure leads to distinct dynamical behavior across the network.  The understanding of the dynamics in heterogeneous networks is in its early stages \cite{Young,HubSync}.

\paragraph{Mutually coupled equations.} Our aim is to study the behavior of solutions of mutually coupled systems with interaction structure given by a heterogeneous random graph.  Consider the set of $n$ nonautonomous linear equations
\begin{equation}\label{unmod1}
\dot{x}_i = V_i(t) x_i, \qquad  \mbox{ for } i = 1, \dots, n,
\end{equation}
where, for $ d \ge 1 $, each $V_i : \mathbb{R} \times  \mathbb{R}^d \rightarrow \mathbb{R}^d$ is a continuous and bounded linear operator. We also assume that these equations are unstable, with nontrivial solutions diverging exponentially fast.

We are interested in the changes in dynamics once the equations are coupled. 
We consider the following one-parameter family of $n$ coupled equations
\begin{equation}\label{eq_acopladas}
\dot{x}_i = V_i(t) x_i + \alpha \sum_{j=1}^n A_{ij} [ H(x_j) - H(x_i)  ], 
\end{equation}
where $\alpha$ is the overall coupling strength, $H$ is a positive-definite matrix, and $A= (A_{ij})$ is the adjacency matrix describing graph connectivity, so that $A_{ij} = 1$ if $i$ receives a connection from $j$ and $A_{ij} = 0$ otherwise.   The degree of the $i^\text{th}$ node to be the number of connections it receives. We will focus on heterogeneous graphs. To be precise, if $\kappa_j < \kappa_i$ denote node degrees in different subnetworks,  then heterogeneity means that  $\Gamma \kappa_i^{1-\gamma} \ge  \kappa_{j}$ for some $0< \gamma <1$ and $\Gamma>0$.  These models have hierarchical organization with modular structures. 

We aim at understanding the dynamics of almost every heterogeneous connection structure $A$.  
The combination of the probabilistic point of view in graphs with the theory of exponential dichotomies  makes it possible to charaterize the dynamics of a large set of networks. 
Our main results show that, for large random graphs, as $\alpha$ is increased, there is a 
bifurcation leading to the  appearance of stable subspaces of solutions. 
Furthermore, the dimension of the stable subspace is determined by the graph structure. 
Loosely speaking, for a suitable range of coupling strength $\alpha$, we have that 
$$
\mbox{ dimension of stable subspace  }  = d \times \ell,  
$$
where $d$ is the dimensional of the solution space of the isolated equation, and $\ell$
is the number of highly connected nodes in the graph. Moreover, if the highly connected nodes are in distinct connectivity regimes, then we prove that there is a cascade of bifurcations increasing the dimension of the stable subspace of solutions.  The precise statements of our results can be found in theorem~\ref{MainThm}
and theorem~\ref{Bifurcativo}.

\section{Preliminaries}

In this section, we provide the basic ingredients for the statement of our results. 

\subsection{Notation} 

We use the small ``$o$'' notation for the asymptotic behavior $n \rightarrow \infty$. 
We write $f(n) = o(1)$ if $f(n)$ goes to zero as $n$ tends to infinity.

We endow the vector space $\mathbb{R}^d$ with the usual Euclidean inner product and the associated Euclidean norm. The space of linear operators is equipped with the induced operator norm. For a continuous family of bounded operators 
$V : \mathbb{R}_+ \times \mathbb{R}^d \to \mathbb{R}^d$, we consider the uniform norm 
$$
\| V \| = \sup_{t>0} \| V(t) \|.
$$
The identity is denoted by $I_d$.

\subsection{Nonautonomous Linear Equations}

We introduce now the concept of exponential dichotomy for a linear differential equation. 
We follow closely \cite{Coppel, Martin}. 

Consider the nonautonomous linear equation
\begin{equation}\label{model}
\dot{x} = V (t) x,
\end{equation}
where $V : \mathbb{R} \times \mathbb{R}^d \rightarrow \mathbb{R}^d$ is a continuous and bounded linear operator. 
We denote by $T(t,s)$ the associated evolution operator, which describes how the solution evolves in time: $ x(t) = T(t,s) x(s) $. 
Clearly,
$$
T(t,t) = I_d \qquad \mbox{  and  } \qquad T(t,s)T(s,r) = T (t,r), \quad \forall \, t,s,r \in \mathbb{R}
$$

The following definition will be central for our study. 

\begin{definition}[Exponential Dichotomy] We say that the linear equation (\ref{model})
admits an exponential dichotomy in the half line $\mathbb{R}_+$ if there is a projector $P: \mathbb{R} \times \mathbb{R}^d \rightarrow \mathbb{R}^d$, with 
$$
P(t) T(t,s) = T(t,s)P(s) \quad \forall \, t\ge s > 0,
$$
for which one may find constants $\eta>0$ and $K>0$ such that, for all $ t \ge s \ge 0 $,
\begin{equation}\label{ExpDic}
\big\| T(t,s)P(s)  \big\| \le K e^{-\eta (t-s)} \quad \mbox{  and  }  \quad \big\| T^{-1}(t,s) (I_d - P)(t) \big\| \le Ke^{-\eta(t-s)}.
\end{equation}
\end{definition}

The exponential dichotomy implies that there is a stable subspace of solutions tending to zero uniformly and exponentially as time goes to infinity. In the complementary subspace,  solutions tend to infinity uniformly and exponentially as time goes to infinity.   

\subsection{Random Graphs}

We will consider random graphs of $n$ nodes modelled by a probability space consisting of the set of labelled undirected graphs of $n$ nodes, endowed with the power-set sigma-algebra, and a probability measure.  We will use a random graph model and terminology from references~\cite{Bollobas,Chung}. This model is an extension of the Erd\"os-R\'enyi model for random graphs with a general degree distribution.  Concerning the terminology, we will adopt the term ``ensemble'' instead of the longer  expression ``probability space''.  

The main point of the model consists in to prescribe the expected values of the node degrees. For convenience, any given sequence of expected degrees $\bm{w}_n = (w_1 , w_2 ,  \cdots, w_n )$ is supposed to verify 
$$ \max_{1\le k\le n} w_k=w_1 \ge w_2 \ge \cdots \ge w_n > 0.
$$
We consider thus an ensemble of random graphs  $G(\bm{w}_n)$ in which an edge between nodes $ i $ and $ j $ is independently assigned 
with success probability 
$$
p_{ij} = \frac{w_i w_j}{\sum_{k=1}^n w_k}.
$$ 
In order to ensure that $p_{ij} \le 1$, we assume that $\bm{w}_n$ is chosen so that  
\begin{equation} \label{DeltaRho}
\big(\max_{1\le k \le n} w_k\big)^2 \le \sum_{k=1}^n w_k.
\end{equation}

A realisation of a graph in the ensemble ${G}(\bm{w}_n)$ is encoded in the adjacency matrix $A = (A_{ij})$ with $(0,1)$-entries determining the connections among nodes of the graph.  The degree $\kappa_i$ of the $i^\text{th}$ node is the number of connections that it receives:  
$$
\kappa_i = \sum_{j=1}^n A_{ij}.
$$ 
Notice that $\kappa_i$ is a random variable whose expected value is exactly the prescribed quantity $w_i$. In particular,  $w_1 = \max_{1\le i\le n} w_i$ is the largest expected value of a degree. 

\paragraph{Network property.}  We say that a property $\mathfrak P$ holds {\it almost surely} if the probability  
that $\mathfrak P$ holds tends to $1$ as $n$ goes to infinity. The assertion {\it almost every graph} in 
${G}(\bm{w}_n)$ has a property $\mathfrak P$ shall be understood as the proportion of all labelled graphs 
of order $n$ that satisfy $\mathfrak P$ tends to 1 as $ n $ goes to infinity.

This asymptotic probabilistic viewpoint naturally motivates us to work with sequences of random graphs and 
thus with sequences of expected degrees. In a rigorous way, we should write $ w_i(n) $ to highlight
the dependence in $ n $ of the expected degree of the $i^\text{th}$ node. In order to avoid a heavy notation,
we leave this dependence implicit. Hence, by imposing additional assumptions on the prescribed expected degrees, 
we focus our attention on the following suitable sequences of heterogeneous graphs.
\begin{definition}[Strong Heterogeneity]
\label{H1} 
For a non decreasing function $ \ell : \mathbb N \to \mathbb N $ and constants $ \theta, \gamma \in (0,1) $, 
we say that a sequence of ensembles $\{{G}(\bm{w}_n)\}_{n\ge1}$ is $(\ell, \theta, \gamma)$-strongly heterogeneous 
if the following hypotheses are satisfied.  
\begin{description}
\item{[H0]} {\it Cardinality of hubs:} there exists a universal constant $ \Gamma_0 > 0 $ such that  
$$
\ell(n)  < \Gamma_0 \big(\max_{1\le k \le n} w_k\big)^\theta.
$$ 
\item{[H1]} {\it Massively connected hubs:} there exists a constant $ c_0 \in (0,1/2] $ such that 
$$
\liminf_{n \to\infty} \, \min_{1 \le i \le \ell(n)} \, \frac{w_{i}}{\max_{1\le k \le n} w_k} =
\liminf_{n \to\infty} \, \frac{w_{\ell(n)}}{\max_{1\le k \le n} w_k}= 2 c_0.
$$
\item{[H2]} {\it Poorly connected nodes:} there exist universal constants $ \Gamma_1, \Gamma_2 > 0 $ and $ \beta > 0 $ such that,
for every $ i \in \{\ell(n) +1, \ldots, n\} $, 
$$
\Gamma_1 \big(\log n\big)^{1+\beta}  <  w_{i}  < \Gamma_2 \big(\max_{1\le k \le n} w_k \big)^{1- \gamma}.
$$
\end{description}
\end{definition}

Notice that $ \ell(n) $ indicates the number of hubs of the graph, that is, of highly connected nodes.  
The parameter $ \theta $ restricts thus their amount.
Besides, thanks to the hypothesis [H2], the constant $ \gamma $ controls the scale separation between low degree nodes and hub nodes.
By abuse of notation, we say that any element ${G}(\bm{w}_n)$ of such a sequence is $(\ell, \theta, \gamma)$-strongly heterogeneous.
We denote this ensemble of heterogeneous random graphs by $G_{\ell, \theta, \gamma}(\bm{w}_n)$. A relevant subclass of heterogeneous
graphs is introduced below.

\begin{definition}[Hubs in Distinct Regimes]
We say that an $(\ell,\theta,\gamma)$-strongly heterogeneous sequence of ensembles $\{G_{\ell,\theta,\gamma}(\bm{w}_n)\}_{n\ge 1}$ has 
hubs in distinct regimes if the additional hypothesis is verified:
\begin{description}
\item{[H1']} there exist sequences of constants $\{\sigma_i\}, \{\tau_i\} \subset (0,1]$ such that both are strictly decreasing and, 
for any fixed index $ i < \ell(n) $,
$$ 
\liminf_{n\to\infty} \frac{w_i}{\max_{1\le k \le n} w_k} > \sigma_i > \tau_i > \limsup_{n\to\infty} \frac{w_{i+1}}{\max_{1\le k \le n} w_k}.
$$ 
\end{description}
\end{definition}

We commit again abuse of notation by extending such a designation to any element of the sequence, which will be denoted by 
$G_{\ell,\theta,\gamma}'(\bm{w}_n)$. 
 
\section{Main Theorems and Discussion}

Consider the uncoupled equations~\eqref{unmod1}. Due to the asymptotic nature of our analysis, we assume from now on that
\begin{equation*}
\max_{i \ge 1} \| V_i \| < + \infty.
\end{equation*} 
The unique solution of each equation can be represented in terms of the transition matrix 
$$
x_i(t) = T_i(t,s) x_i(s), \qquad i = 1, \ldots, n.
$$
When we say that solutions are unstable, we mean that there are constants $\eta_0>0$ and $K_0>0$ such that, for all $ t > s $,
\begin{equation}\label{mesmasconstantes}
\| T_i^{-1}(t,s) \| \le K_0 \, e^{- \eta_0 (t-s)}, \qquad  i = 1, \ldots, n.
\end{equation}
Notice that we suppose all evolution operators share the same constants $ \eta_0 $ and $ K_0 $.

To state our results, it will be convenient to represent coupled equations in a single block form. 
Consider 
$$
X =  \mbox{col}(x_1, \cdots, x_n),
$$
where col denotes the vectorization formed by stacking the column vectors $x_i$ into a
single column vector. Moreover, we denote 
$$
V(t) = \bigoplus_{i=1}^n V_i(t) = \mbox{diag}(V_1(t), \dots, V_n(t)).
$$
Then the coupled equations \eqref{eq_acopladas} can be recast into a block form
\begin{equation}\label{md1}
\dot X = [V(t) - \alpha L \otimes H] X,
\end{equation}
where $L = (L_{ij}) $ is a combinatorial laplacian given by $L_{ij} = \delta_{ij} \kappa_i - A_{ij}$
(as usual, $\delta_{ij}$ stands for the Kronecker delta),
and $\otimes$ is the Kronecker product \cite{Lancaster}. The unique solution of equation~\eqref{md1} can be 
represented in terms of the transition matrix 
$$
X(t) = T(t,s) X(s).
$$

For $\alpha=0$, the  equations are uncoupled and have only unstable solutions. 
Our main results show that stable solutions appear when these equations are coupled in heterogeneous random graphs and that
increasing the coupling strength $\alpha$ leads to a sequence of bifurcations.

\begin{theorem}\label{MainThm}
Consider the ensemble $G_{\ell, \theta, \gamma}(\bm{w}_n)$, with $\theta < (3 - \sqrt{5})/2 $ and $\gamma  > (\sqrt{5} - 1)/2$. 
Then, there are constants $ C, c >0 $ and an integer $ N_0> 0$ such that,  for all $ n> N_0 $, whenever
$$
c < \alpha \max_{1\le k \le n} w_k < C \big( \log n \big)^\gamma,
$$ 
with probability at least $ 1 - n^{-1/2} - 2n^{-1/5} $, the coupled equations \eqref{md1} on a graph of $G_{\ell, \theta, \gamma}(\bm{w}_n)$ admit an exponential dichotomy:
for positive constants $K$ and $\eta$ and a projector $P$ that commutes with $T$, for all $t\ge s\ge 0$,
\begin{equation*}
\big\| T(t,s)P(s)  \big\| \le K e^{-\eta (t-s)} \quad \mbox{  and  }  \quad \big\| T^{-1}(t,s) (I_d - P)(t) \big\| \le Ke^{-\eta(t-s)}.
\end{equation*}
Moreover, in such a situation, the dimension of the stable subspace is determined by the network structure
$$
\mbox{\rm rank } P = d \times  \ell(n).
$$
\end{theorem}

Roughly speaking, the constants $c$ and $C$ in our theorem \ref{MainThm} are given by two distinct mechanisms. The constant $c$ comes from the fact that we wish to guarantee the existence of the stable subspace whereas the constant $C$ comes from fact that the complementary subspace of unstable solution must have a uniform exponent divergence. And thereby we ensure the existence of the dichotomy.

On a heterogeneous random graph a natural coupling parameter is given by 
$$ 
\alpha = \frac{\alpha_0}{\max_{1\le k \le n} w_k}. 
$$ 
We regard the parameter $\alpha_0$ as the normalized coupling strength capturing the dynamics at the highly connected nodes.
In this case, the previous theorem can be restated in the following form.

\begin{theorem}
For the ensemble $G_{\ell, \theta, \gamma}(\bm{w}_n)$, with $\theta < (3 - \sqrt{5})/2 $ and $\gamma  > (\sqrt{5} - 1)/2$, there exists a positive constant 
$ c = c(H,\max_{i\ge 1} \|V_i\|, c_0)$ 
such that if 
$$ 
\alpha_0 > c 
$$
then the coupled equations~\eqref{md1} on almost every graph of $G_{\ell, \theta, \gamma}(\bm{w}_n)$ admit an exponential dichotomy,
in which the dimension of the stable subspace is exactly $ d \times \ell(n) $.
\end{theorem}

For the case of hubs in distinct regimes, we highlight a bifurcation-type result.

\begin{theorem}\label{Bifurcativo}
Given an ensemble $G_{\ell,\theta,\gamma}'(\bm{w}_n)$, with $\theta < (3 - \sqrt{5})/2 $ and $\gamma  > (\sqrt{5} - 1)/2$,
there are constants $\bar{C}, \bar{c} > 0$ and an integer $ N_0 > 0 $ such that for any $ n > N_0 $, if
$$ \frac{\bar c}{\sigma_i} < \alpha \max_{1\le k \le n} w_k < \frac{\bar C}{\tau_i} 
\qquad \text{for some index } \, i < \ell(n), $$
then, with probability greater or equal to $ 1 - n^{-1/2} - 2n^{-1/5} $, the coupled equations~\eqref{md1} on a graph of  
$G_{\ell,\theta,\gamma}'(\bm{w}_n)$ admit an exponential dichotomy, in which the dimension of the stable space is $ d \times i $.
\end{theorem}

The constants $ \bar C = \bar C(H, \eta_0) $ and $ \bar c = \bar c(H, \max_{i \ge 1} \| V_i \|) $ may be explicitly determined (see section 5).
If $ \bar c / \bar C < \sigma_i / \tau_i $ for all $ i < \ell(n) $, notice that, as $ \alpha $ is increased, the 
system exhibits a cascade of bifurcations, characterized by the increasing of the dimension of the stable subspace of solutions.  
For an illustration, suppose that the elements of the problem are chosen so that $ \bar c < \bar C $. In this case, we may assume
$ \sigma_i = \tau_i $. Note that $ \sigma_{i-1} $ and $ \sigma_i $ control then the proportion of the hub $ i $ has with respect to
the main hub: 
$$ \sigma_i \max_{1\le k \le n} w_k < w_i < \sigma_{i-1} \max_{1\le k \le n} w_k $$ 
for $ n $ large enough. We assume 
in addition that $ \bar c / \bar C = \sigma_i / \sigma_{i-1} $ for all $ i $. Thus, for the coupling constant 
$ \alpha = \alpha_0 / \max_{1\le k \le n} w_k $, the inequalities
$$ \bar c < \alpha_0 \sigma_i < \bar C $$
imply that, on almost every graph of $ G_{\ell,\theta,\gamma}'(\bm{w}_n) $, the coupled equations~\eqref{md1} admit an exponential 
dichotomy, in which the dimension of the stable space is $ d \times i $. Besides, as $ \alpha_0 $ is increased, 
the global bifurcation occurs from the transition of an interval control condition 
$$ 
\alpha_0 \sigma_{i+1} < \bar c < \alpha_0 \sigma_i < \bar C < \alpha_0 \sigma_{i-1} 
$$
to the next one  
$$
\alpha_0 \sigma_{i+2} < \bar c < \alpha_0 \sigma_{i+1} < \bar C < \alpha_0 \sigma_i. 
$$

We conclude this section describing the main ideas of the proof of theorem~\ref{MainThm}. Its proof will be given in section 4.
In section 5, we point out which minor changes have to be made in order to prove theorem~\ref{Bifurcativo}.

\paragraph{Strategy of the proof of theorem~\ref{MainThm}.} We rewrite the block form~\eqref{md1} of the coupled equations as
$$
\dot{X} = [V(t) - \alpha D \otimes H ] X +\alpha  (A \otimes H) X,  
$$
where $D = $ diag$(\kappa_1, \dots, \kappa_n)$ is the matrix of degrees. 
Our strategy is to obtain the existence of dichotomies by persistence arguments. Essentially,
the proof consists of three steps: 
\begin{itemize}
\item[\emph{i})] Notice that the ensemble of random graphs has concentration properties: the actual degrees $\kappa_i$'s are almost surely described by the expected degrees $w_i$'s. 
\item[\emph{ii})] Treat $\alpha  ( A \otimes H) $ as a perturbation. Using the block form of 
$$
\dot{Y} = 
\left(
\begin{array}{ccc}
V_1(t) - \alpha \kappa_1  H & \cdots & 0 \\
\vdots & \ddots & \vdots \\
0 & \cdots & V_n(t) - \alpha \kappa_n H
\end{array}
\right) Y
$$ 
together with the concentration properties of the degrees and scale separation [H2] of the ensemble 
$G_{\ell, \theta, \gamma}(\bm{w}_n)$, 
we prove that, for a suitable coupling strength $\alpha$, the first $\ell(n)$ 
blocks associated with the highly connected nodes are exponentially 
stable whereas the remaining ones are unstable.
\item[\emph{iii})] Include the term $\alpha (A \otimes H) $ and use persistence of dichotomies. 
The challenge here is to proof that this coupling term is small. We use the concentration properties of $G_{\ell, \theta, \gamma}(\bm{w}_n)$ and the additional conditions on the parameters scale separation $\gamma$ and cardinality of hubs $\theta$ to show that, in the limit of large graphs, $\alpha \| A \otimes H \|$ can be made arbitrary small with respect to the dichotomy parameters. Then, we apply the persistence of exponential dichotomies to obtain the persistence of the stable subspace 
of solutions.  
\end{itemize}

\section{Proof of Theorem~\ref{MainThm}}

Before providing the details of the proof, we need some auxiliary results. We group them according to the research domain.

\subsection{Random Graphs}

We start estimating the actual degrees with respect to the expected degrees. 
The next result will be very useful in such an analysis. For a proof, see the demonstration of lemma~5.7 in~\cite{Chung}.

\begin{propo}\label{kiwi}
Graphs in an ensemble $G(\bm{w}_n)$ have degree concentration property in the sense that:
$$ \mathbb P\Big[|\kappa_i - w_i| \le 2\sqrt{\log n}\sqrt{\max\{w_i, \log n\}} \quad \forall \, i= 1,\ldots, n \Big] \ge 1- 2n^{-1/5}. $$
In particular, if $ \log n / \min_{1\le k \le n} w_k $ tends to zero as $ n $ goes to infinity, then
$$
\kappa_i = w_i \left[ 1 + o(1) \right], \quad \forall \, i = 1, \ldots, n,
$$
holds almost surely.
\end{propo}

The previous degree concentration property allows us to highlight interesting facts. 

\begin{coro}\label{concentracao}
For graphs in a strongly heterogeneous ensemble $ G_{\ell, \theta, \gamma}(\bm{w}_n) $, whenever $ n $ is sufficiently large,
the probability of the event
\begin{equation*}
\max_{\ell(n) < i\le n}\kappa_i<\frac{3}{2}\Gamma_2\big(\max_{1\le k\le n}w_k\big)^{1-\gamma} \quad \text{and} \quad 
\min_{1\le i \le\ell(n)}\kappa_i>\frac{1}{2}c_0\max_{1\le k\le n}w_k  
\end{equation*}
is at least $ 1 - 2 n^{-1/5} $. 
\end{coro}

\begin{proof}
Thanks to the previous proposition, from hypothesis [H2], for $ n $ sufficiently large, the probability of the event 
$$ \kappa_i \le w_i\big( 1 + 2\sqrt{\frac{\log n}{w_i}} \big) 
< \Gamma_2\big(\max_{1\le k\le n}w_k\big)^{1-\gamma} \big(1 + \frac{2}{\sqrt{\Gamma_1 (\log n)^\beta}} \big), 
\quad \forall \, i = \ell(n) + 1, \ldots, n, $$
is greater or equal to $ 1 - 2 n^{-1/5} $.
Moreover, hypothesis [H1] guarantees that, for $ n $ large enough, the event 
$$ \kappa_i \ge w_i \big( 1 - 2\sqrt{\frac{\log n}{w_i}}\big) > c_0 \max_{1\le k\le n}w_k \, \big(1-\frac{2}{\sqrt{\Gamma_1(\log n)^\beta}}\big), 
\quad \forall \, i = 1, \ldots, \ell(n), $$
occurs simultaneously with at least the same probability.
\end{proof}

The proof of the next corollary is similar and will be omitted.

\begin{coro}\label{concentracaobifurcada}
For graphs in a strongly heterogeneous ensemble with hubs in distinct regimes $G_{\ell, \theta, \gamma}'(\bm{w}_n)$, if $ n $ is large
enough, then the event 
\begin{equation*}
\max_{j < i\le n}\kappa_i<\frac{3}{2} \tau_j \max_{1\le k\le n}w_k \quad \text{and} \quad 
\min_{1\le i \le j}\kappa_i>\frac{1}{2} \sigma_j \max_{1\le k\le n}w_k, \qquad j < \ell(n),   
\end{equation*}
occurs with probability at least $ 1 - 2 n^{-1/5} $.
\end{coro}

Given a sequence of expected degrees $\bm{w}_n = (w_1, \ldots, w_n)$, the second-order average degree is given by
$$
\Delta := \frac{\sum_{k=1}^n w_k^2}{\sum_{k=1}^n w_k}.
$$
This constant plays an important role for the characterisation of the ensemble $G(\bm{w}_n)$.

\begin{propo} \label{CLV}
Suppose that the largest expected degree satisfies $ \max_{1\le k \le n} w_k \ge \log n$.
Let $ \lambda_\text{max} $ denote the largest eigenvalue of the adjacency matrix associated with a random graph in $G(\bm{w}_n)$. 
Then, the probability of the event
$$ \lambda_\text{max} \le \Delta + \frac{3}{2}\sqrt{\log n \, \max_{1\le k \le n} w_k} +
\sqrt{\frac{1}{4} \log n \, \max_{1\le k \le n} w_k + 3(\Delta + \log n)\sqrt{\log n \, \max_{1\le k \le n} w_k}} $$
is greater or equal to $ 1 - n^{-1/2} $. 
\end{propo}
\begin{proof}
This result is actually a minor modification of lemma~3.2 in~\cite{ChungLuVu2003}.
Just repeat the same proof there using 
$$ a = \frac{1}{2} \sqrt{\log n \, \max_{1\le k \le n} w_k}  + \sqrt{\frac{1}{4} \log n \, \max_{1\le k \le n} w_k + 3(\Delta + \log n)\sqrt{\log n \, \max_{1\le k \le n} w_k}}. $$
\end{proof}

The previous proposition leads to an important result on the control of the normalized perturbation size for heterogeneous graphs. 
The following statement will be of fundamental importance for us. 

\begin{propo}\label{m}
Consider a strongly heterogeneous ensemble $G_{\ell, \theta, \gamma}(\bm{w}_n)$ with $\theta < \delta/2 < (3-\sqrt{5})/2 $ and 
$\gamma > 1 - \delta/2 > (\sqrt{5}-1)/2 $, where $ \delta \ge 3/4 $.  
Then, for $n$ large enough, with probability at least $ 1 - n^{-1/2} $, the largest eigenvalue $ \lambda_\text{max} $ of the adjacency matrix associated with a random graph in $G_{\ell, \theta, \gamma}(\bm{w}_n)$ verifies
$$
\lambda_{\text{max}} \le  5 \big(\max_{1\le k \le n} w_k\big)^\delta.
$$
\end{propo}

\begin{proof}
Assuming strong heterogeneity as above, we first show that 
\begin{equation}\label{majorante1}
 \Delta < \big(\max_{1\le k \le n} w_k\big)^\delta 
\end{equation}
for $ n $ large enough. Notice that, from~\eqref{DeltaRho} and hypotheses [H0] and [H2], we have 
$$ \big(\max_{1\le k \le n} w_k\big)^2 < \Gamma_0 \big(\max_{1\le k \le n} w_k\big)^{1 + \theta} + 
\Gamma_2 \, n \big(\max_{1\le k \le n} w_k\big)^{1-\gamma}. $$
Since $ \theta < 1 $, for $ n $ large enough, we obtain that
$$  \big(\max_{1\le k \le n} w_k\big)^2 < 2 \Gamma_2 \, n \big(\max_{1\le k \le n} w_k\big)^{1-\gamma}. $$
Using again hypotheses [H0] and [H2], for $n$ sufficiently large, we get
$$
 \sum_{k=1}^n w_k > n \qquad \text{and} \qquad
\sum_{k=1}^n w_k^2 < \Gamma_0 \big(\max_{1\le k \le n} w_k\big)^{2 + \theta} + 
\Gamma_2^2 \, n \big(\max_{1\le k \le n} w_k\big)^{2-2\gamma},
$$
so that
$$ \Delta < \Gamma_0 \, \frac{1}{n} \big(\max_{1\le k \le n} w_k\big)^{2 + \theta}
+ \Gamma_2^2 \big(\max_{1\le k \le n} w_k\big)^{2-2\gamma}. $$
Hence, for $ n $ large enough,
$$ \Delta < 2\Gamma_0\Gamma_2 \big(\max_{1\le k \le n} w_k\big)^{1+\theta-\gamma}
+ \Gamma_2^2 \big(\max_{1\le k \le n} w_k\big)^{2-2\gamma} $$
and the desired inequality follows from the choice of the parameters $ \theta $ and $ \gamma $.

Note now that, when $ \gamma \ge 1/2 $, clearly  
$ \log n < (\Gamma_1/\Gamma_2)^{\frac{1}{2(1-\gamma)}} \big(\log n\big)^{\frac{1+\beta}{2(1-\gamma)}} $
for $n$ large enough. Thus, from hypothesis [H2], it follows that
\begin{equation}\label{majorante2}
\log n < \sqrt{\max_{1\le k \le n} w_k}
\end{equation} 
for $n$ sufficiently large.

From inequalities~\eqref{majorante1} and~\eqref{majorante2}, and from the fact that $ \delta \ge 3/4 $, 
it is easy to deduce 
for $n$ large enough 
\begin{multline*}
 \Delta + \frac{3}{2}\sqrt{\log n \, \max_{1\le k \le n} w_k} +
\sqrt{\frac{1}{4} \log n \, \max_{1\le k \le n} w_k + 3(\Delta + \log n)\sqrt{\log n \, \max_{1\le k \le n} w_k}} \le \\
\le 5 \big(\max_{1\le k \le n} w_k\big)^\delta. 
\end{multline*}
Therefore, the result follows from proposition~\ref{CLV}.
\end{proof}

\subsection{Exponential Dichotomies}

One of the most important properties of exponential dichotomies is their roughness. 
In clear terms, they are not destroyed by small perturbations on the matrix entries. 
Our proof relies on this persistence. Therefore, for completeness we state the following result.
A proof can be found in~\cite{Coppel}.

\begin{lemma} \label{roughness}
Suppose that the linear differential equation~\eqref{model} admits an exponential dichotomy~\eqref{ExpDic} on $\mathbb{R}_+$. 
If 
$$
\delta := \sup_{t \in \mathbb{R}_+} \|B(t) \| < \frac{\eta}{4 K^2}, 
$$
then the perturbed equation 
$$
\dot{y} = [V(t) + B(t) ]{y}
$$
has a similar exponential dichotomy in the sense that: there exist a constant $\hat{K}>0$ and a projector $ \hat P $, which preserves the rank
of the original projector $ P $ and commutes with the evolution operator $\hat T$ associated with the perturbed equation, 
such that, for $ t \ge s \ge 0$, 
\begin{equation*}
\big\| \hat{T}(t,s) \hat P(s) \big\| \le \hat{K} e^{-(\eta - 2 K \delta) (t-s)} \quad \text{and} \quad 
\big \| \hat T^{-1}(t,s) (I_d - \hat P) (t) \big \|  \le   \hat{K} e^{-(\eta - 2 K \delta) (t-s)}.
\end{equation*}
\end{lemma}

The next propositions will be important in our proof.

\begin{propo}\label{Stab}
Consider the equation 
\begin{equation}\label{p1}
\dot{y} = [V(t) - \alpha H]  y,
\end{equation}
where $H$ is a positive-definite matrix. 
Let $ \hat T $ be the associated evolution operator. 
Then, there exist a constant $\hat K_H >0$ that only depends on $ H $ such that
$$
\| \hat T(t,s) \| \le \hat K_H e^{-(\alpha \lambda_H - \hat K_H \| V \|)(t-s)}, \qquad \forall \, t \ge s,
$$
where $\lambda_H > 0$ is the smallest eigenvalue of $H$.
\end{propo}

\begin{proof}  First we solve $ \dot{x} = - \alpha H  x$. 
Clearly,  in this case $T(t,s) = e^{-\alpha (t-s) H}$ is the associated evolution operator. 
Notice that there is a positive constant $ \hat K_H $ such that $\| T(t,s) \| \le \hat K_H \exp[ - \alpha \lambda_H(t-s) ]$
for all $ t \ge s $, where $\lambda_H > 0$ is the smallest eigenvalue of $H$. 
For the full equation~\eqref{p1}, the variation-of-constants formula yields
$$ \hat T(t,s) - T(t,s) = \int_s^t T(t, \tau) V(\tau) \hat T(\tau, s) \, d\tau. $$
Hence, we obtain
$$
\| \hat T(t,s) \| \le \hat K_H e^{  - \alpha \lambda_H (t-s) }  + 
\hat K_H  \int_s^t e^{ - \alpha \lambda_H (t - \tau)} \,  \| V(\tau) \| \, \| \hat T(\tau, s) \| \, d\tau.
$$
Now introducing $\omega_{s,t}(u) = \omega(u) = \| \hat T (u,s)  \| \exp[- \alpha \lambda_H (t-u) ] $, we have 
$$
\omega (t) \le \hat K_H \omega (s) + \hat K_H \| V \| \int_s^t \omega(\tau) \, d\tau.  
$$
Using a Gronwall estimate, we conclude that 
$$
\| \hat T(t,s) \| \le \hat K_H e^{-( \alpha \lambda_H - \hat K_H \| V \| ) (t-s)} .
$$
\end{proof}

\begin{propo}\label{Unst}
For the equation $\dot{x} = V(t)  x$, suppose that the associated evolution operator $ T $ verifies
$\| T^{-1}(t,s) \| \le K e^{-\eta (t-s)}$ for all $ t\ge s$. Then, the perturbed equation
$$ \dot{y} = [V(t) + H]  y $$ 
admits an evolution operator $\hat{T}$ satisfying 
$$ \| \hat{T}^{-1}(t,s) \| \le K e^{-(\eta - K \|H\| ) (t-s)}, \qquad \forall \, t \ge s. $$
\end{propo}

\begin{proof} 
Notice that the respective evolutions operators verify the following partial differential equations
$$ \partial_t T^{-1}(t,s) = - T^{-1}(t,s) V(t) \qquad \text{and} \qquad 
\partial_t \hat T^{-1}(t,s) = - \hat T^{-1}(t,s) [V(t) + H]. $$
Using variation of constants, we then obtain
$$ \hat T^{-1}(t,s) - T^{-1}(t,s) = \int_s^t \hat T^{-1}(\tau, s) [-H] T^{-1}(t,\tau) \, d\tau. $$
Hence, by the triangle inequality, we have
$$
\| \widehat T^{-1}(t,s) \| \le K e^{-\eta(t-s)} + K \| H \| \int_s^t \| \widehat T^{-1}(\tau,s)\|   e^{-\eta(t-\tau)} \, d\tau.
$$
Following the same steps as in the previous proposition, now with $ \omega_{s,t}(u) = \omega(u) = \| \hat T^{-1}(u,s) \| \exp[-\eta (t - u)] $, 
we obtain the result. 
\end{proof}

\subsection{The Proof}

We follow the strategy presented at the end of section 3.
 
\begin{proof}[Proof of Theorem~\ref{MainThm}] \

\paragraph{Step i. --} In corollary~\ref{concentracao}, we have already established in a precise way how the prescribed expected degrees $w_i$'s 
almost surely determine actual degrees $\kappa_i$'s.

\medskip

\paragraph{Step ii. --}
We can rewrite the equation~\eqref{md1} as follows
\begin{equation}
\dot{X} =   [ \Omega(t,\alpha) + \alpha A \otimes H ] X,
\label{block}
\end{equation}
where $ {\Omega}(t,\alpha) := \bigoplus_{i=1}^n \left[ V_i(t) - \alpha \kappa_i H \right] $. 
Now consider the system 
\begin{equation}\label{acopladasimples}
\dot{Y} =   {\Omega}(t,\alpha) Y.
\end{equation}
Since the system is block diagonal, we can solve each block independently. 
Let 
$$
\hat{T}(t,s) = \bigoplus_{i=1}^n \hat{T}_i(t,s)
$$
 be the associated evolution operator. 
On the one hand, applying proposition~\ref{Stab} to the first $ \ell(n) $ blocks, we conclude that  
$$
\|\hat{T}_i(t,s)\|\le\hat K_H\exp\big[-\big(\alpha\lambda_H\min_{1\le i\le\ell(n)}\kappa_i-\hat K_H\max_{1\le i \le\ell(n)}\| V_i \|\big)(t-s)\big], 
\qquad \forall \, 1\le i \le \ell(n).
$$
Recall that by hypothesis $ \max_{i \ge 1}\| V_i \| < +\infty $.
On the other hand, proposition~\ref{Unst} shows that the remaining blocks will verify 
$$
\| \hat{T}^{-1}_i(t,s) \| \le K_0 \exp\big[- \big(\eta_0 - \alpha \| H \|\max_{\ell(n) < i \le n} \kappa_i \big) (t-s)\big], \qquad  \forall \,  
\ell(n) <  i \le n,
$$
where $ \eta_0 $ and $ K_0 $ are the universal constants introduced in~\eqref{mesmasconstantes}.
These observations lead us to consider the operator $\hat{T}(t, s)$ associated with equation~\eqref{acopladasimples} 
in a block form with respect to the direct sum $\mathbb{R}^{d \times \ell(n)} \oplus \mathbb{R}^{d \times (n - \ell(n)) }$, namely 
$$
\hat{T} (t, s) =
\left(
\begin{array}{cc}
S (t, s) & 0 \\ 
0 & U (t, s)
\end{array}
\right),
$$ 
as well as the natural projectors
$$
\hat{P} =
\left(
\begin{array}{cc}
I_d & 0 \\ 
0 & 0
\end{array}
\right)
\qquad \text{ and } \qquad
Id - \hat P = 
\left(
\begin{array}{cc}
0 & 0 \\ 
0 & I_d
\end{array}
\right).
$$ 
Since the operator norm induced by the Euclidean norm has the property 
$$
\|\hat T(t,s)\| = \max_{1 \le i \le n}  \|\hat T_i(t,s)\| \quad \text{and} \quad
\|\hat T^{-1}(t,s)\| = \max_{1 \le i \le n}  \|\hat T^{-1}_i(t,s)\|, 
$$
in order to characterize an exponential dichotomy, the contraction rates for stable and unstable directions must satisfy
\begin{equation}\label{eventotaxas}
\alpha \min_{1 \le i \le \ell(n)} \kappa_i  > \frac{\hat K_H \max_{i\ge1} \| V_i \|}{\lambda_H}  
\qquad \mbox{ and } \qquad 
\alpha \max_{\ell(n) < i \le n} \kappa_i < \frac{\eta_0}{\| H \|}. 
\end{equation}
Suppose from now on that $ n $ is large enough so that
$$ \big(\max_{1\le k\le n} w_k \big)^\gamma > \big(\log n\big)^\gamma > 
\frac{12 \Gamma_2 \| H \| \hat K_H \max_{i \ge 1} \|V_i\|}{c_0 \eta_0 \lambda_H}. $$
Thus, defining 
\begin{equation*}
c:= \frac{4 \hat K_H \max_{i\ge 1} \|V_i\|}{c_0 \lambda_H}  \qquad \mbox{  and  }  \qquad C := \frac{\eta_0}{3 \Gamma_2 \| H \|}, 
\end{equation*}
let $ \alpha > 0 $
be such that
$$ c < \alpha \max_{1\le k\le n} w_k < C\big(\log n\big)^\gamma. $$
Notice now that
\begin{align*}
\frac{c}{2} < \alpha \max_{1\le k\le n} w_k \quad & \Rightarrow \quad 
\frac{1}{2} \alpha c_0 \max_{1\le k\le n} w_k > \frac{\hat K_H \max_{i\ge 1}\|V_i\|}{\lambda_H}, \\
\alpha \max_{1\le k\le n} w_k < 2C\big(\max_{1\le k\le n} w_k \big)^\gamma \quad & \Rightarrow \quad
 \frac{3}{2} \alpha \Gamma_2 \big(\max_{1\le k\le n} w_k \big)^{1-\gamma} < \frac{\eta_0}{\|H\|}.
\end{align*}
Therefore, from corollary~\ref{concentracao}, for $ n $ sufficiently large, we conclude that the event~\eqref{eventotaxas} occurs with
probability at least $ 1 - 2 n^{-1/5} $. With the same estimate for
the probability, equation~\eqref{acopladasimples} admits thus an exponential dichotomy with constants
\begin{align*}
\hat K & = \max\{ \hat K_H, K_0\} \qquad \text{and} \\ 
\hat \eta & = \min\Big\{\frac{1}{2}  \lambda_H c_0 \alpha \max_{1\le k\le n} w_k - \hat K_H \max_{i\ge 1}\|V_i\|, \,
\eta_0 - \frac{3}{2}  \|H\| \Gamma_2 \alpha \big(\max_{1\le k\le n} w_k \big)^{1-\gamma}\Big\}.
\end{align*}

\medskip

\paragraph{Step iii. --} Now we wish to incorporate back the perturbation $\alpha A\otimes H$. 
Notice that $ \| \alpha A \otimes H \| \le \alpha \| A \| \| H \| $. Moreover, since $ A $ is a real symmetric matrix,
$ \| A \| = \lambda_{\text{max}} $. Thus, in order to apply lemma~\ref{roughness},
we need to estimate the probability of the event $ \alpha \lambda_{\text{max}} < \Lambda $, where
\begin{align*} 
\Lambda & := \frac{1}{4 \hat K^2} \min\Big\{\frac{1}{2}\frac{\lambda_H}{\|H\|} c_0 \alpha \max_{1\le k\le n} w_k - 
\frac{\hat K_H}{\|H\|} \max_{i\ge 1}\|V_i\|, \, \frac{\eta_0}{\|H\|} - 
\frac{3}{2} \Gamma_2 \alpha \big(\max_{1\le k\le n} w_k \big)^{1-\gamma}\Big\} \\
& > \frac{1}{4 \hat K^2} \min\Big\{\frac{\hat K_H}{\|H\|} \max_{i\ge 1}\|V_i\|, \, \frac{1}{2} \frac{\eta_0}{\|H\|} \Big\} > 0.
\end{align*}
Thanks to proposition~\ref{m}, with probability at least $ 1 - n^{-1/2} $, 
$$ \alpha \lambda_{\text{max}} \le 5 \alpha \big(\max_{1\le k \le n} w_k\big)^\delta, \qquad 
\text{where } \, \frac{3}{4} \le \delta < 3 - \sqrt{5}. $$
Using the definition of $ \alpha $ and hypothesis [H2], notice then that
$$
5 \alpha \big(\max_{1\le k \le n} w_k\big)^\delta 
< 5 C \frac{\Gamma_2^\gamma}{\Gamma_1^\gamma}\big(\max_{1\le k \le n} w_k\big)^{\gamma - \gamma^2 + \delta - 1}.
$$ 
By hypothesis $ \gamma > 1 - \delta/2 $. Hence, 
$$
\big(\max_{1\le k \le n} w_k\big)^{\gamma - \gamma^2 + \delta - 1} \le \big(\max_{1\le k \le n} w_k\big)^{-\frac{\delta^2 -6\delta +4}{4}}
$$
tends to zero as $n$ goes to infinity, which concludes the proof.
\end{proof}

\section{Comments on the Proof of Theorem~\ref{Bifurcativo}}

The proof of theorem~\ref{Bifurcativo} follows the same lines of the previous one.
Corollary~\ref{concentracaobifurcada} provides us now the first step.
Moreover, the last step is exactly as before, only with a convenient $ \hat \eta $ as described below.

Concerning the second step, the set of arguments remains unchanged.
However, one first introduces constants
$$ \bar{c} := \frac{3 \hat K_H  \max_{i\ge 1} \|V_i\|}{\lambda_H} \qquad \mbox{  and  } \qquad \bar{C} := \frac{\eta_0}{2 \|H\|}. $$
Note that, if $ \bar{c}/\bar{C} < \sigma_j/\tau_j $ for some $ j < \ell(n) $, then, for all $ n $,  it will be possible to consider
a suitable parameter $ \alpha > 0 $ as in the statement. Applying proposition~\ref{Stab} to the first $ j $ blocks
and proposition~\ref{Unst}
to the remaining ones, one has thus to estimate the probability of the event
\begin{equation}\label{eventotaxasbifurcado}
\alpha \min_{1\le i\le j} \kappa_i > \frac{\hat K_H \max_{i \ge 1} \|V_i\|}{\lambda_H} \qquad \text{and} \qquad 
\alpha \max_{j < i \le n} \kappa_i < \frac{\eta_0}{\|H\|}.
\end{equation}
Corollary~\ref{concentracaobifurcada} and the fact that
$$ \frac{3}{2} \alpha \tau_j \max_{1\le k \le n} w_k < 2 \bar{C} \qquad \text{and} \qquad 
\frac{1}{2} \alpha \sigma_j \max_{1\le k \le n} w_k > \frac{\bar{c}}{3} $$
show that~\eqref{eventotaxasbifurcado} occurs with probability at least $ 1 - 2 n^{-1/5} $. Hence, exponential dichotomy
is found with this estimated probability, being now
$$ \hat \eta 
= \min \Big\{ \frac{1}{2} \alpha \sigma_j \max_{1\le k \le n} w_k - \frac{\bar{c}}{3}, 2\bar{C} - \frac{3}{2} \alpha \tau_j \max_{1\le k \le n} w_k \Big\}
> \min \Big\{\frac{\bar{c}}{6}, \frac{\bar{C}}{2}\Big\} > 0. $$

\end{document}